\documentclass[12pt]{article}
\usepackage{amsmath, amsfonts, amsthm}
\usepackage{amssymb,mathrsfs, verbatim}
\usepackage{graphicx,amstext}
\usepackage[mathcal]{euscript}
\usepackage{mathtools}

\newcommand{\clg}[1]{{\mathcal{#1}}}

\newcommand{\PM}{\mathbb P}
\newcommand{\R}{\mathbb R}
\newcommand{\Z}{\mathbb Z}
\newcommand{\N}{\mathbb N}
\newcommand{\C}{\mathbb C}

\newcommand{\ve}{\varepsilon}
\newcommand{\vp}{\varphi}

\newcommand \loc    {\text{loc}}

\newcommand{\Medint}{\mkern12mu\mbox{\vrule height4pt
         depth-3.2pt
          width5pt}\mkern-16.5mu\int\nolimits}

\newcommand{\tobo}[1]{\raisebox{-1.7ex}{$\stackrel{\textstyle {\longrightarrow}}
                                                    {\scriptstyle {#1}}$}}

\newtheorem{theorem}{Theorem}[section]

\newtheorem{remark}[theorem]{Remark}
\newtheorem{lemma}[theorem]{Lemma}
\newtheorem{corollary}[theorem]{Corollary}
\newtheorem{definition}[theorem]{Definition}
\newtheorem{example}[theorem]{Example}

\begin{document}
\title{On a Characterization of the Rellich–Kondrachov Theorem on Groups} 

\author{Vernny Ccajma$^1$, Wladimir Neves$^{1}$, Jean Silva$^2$}

\date{}

\maketitle

\footnotetext[1]{ Instituto de Matem\'atica, Universidade Federal
do Rio de Janeiro, C.P. 68530, Cidade Universit\'aria 21945-970,
Rio de Janeiro, Brazil. E-mail: {\sl
v.ccajma@gmail.com, wladimir@im.ufrj.br.}}
\footnotetext[2]{ Departamento de Matem\'atica, Universidade Federal
de Minas Gerais. E-mail: {\sl jcsmat@mat.ufmg.br}}


%
%
\begin{abstract} Motivated by an eigenvalue-eigenfunction problem posed in $\R^n \times \Omega$, where $\Omega$
is a probability space, we are concerned in this paper with the Sobolev space on groups. 
Hence it is established an equivalence between locally compact Abelian groups
and the space of solutions to the associated variational problem. Then, we 
study some conditions which characterize in a precisely manner 
the Rellich-Kondrachov Theorem, the principal ingredient to solve
the variational problem.
\end{abstract}

\maketitle



\section{Introduction}
\label{Intro}

The main motivation in this paper to study Sobolev spaces on groups, besides
being an elegant and modern mathematical theory, is the Rellich–Kondrachov Theorem. 
More precisely, we are interested on a characterization of it, which is related to the 
solution of the following eigenvalue-eigenfunction problem: 
Find $\lambda(\theta) \in \mathbb{R}$ and $\Psi(\theta)$ a complex-value function 
satisfying for any $\theta \in \R^n$ fixed, 
\begin{equation}
\label{92347828454trfhfd4rfghjls}
\left\{
\begin{array}{l}
L^\Phi(\theta) [\Psi(z,\omega)]= \lambda \ \Psi(z,\omega), \hspace{40pt} \text{in $\R^n \times \Omega$}, 
\\[5pt]
\hspace{32pt} \Psi(z, \omega) = \psi {\left( \Phi^{-1} (z, \omega), \omega \right)}, \quad \text{$\psi$ is a stationary function},
			\end{array}
			\right.
\end{equation}
where $(\Omega, \mathbb{P})$ is a probability space, and the linear operator $L^\Phi(\theta)$ is defined by
\begin{equation}
\label{EqEsp}
\begin{aligned}
L^\Phi(\theta)[\cdot]:=& -\big( {\rm div}_{\! z} + 2i\pi \theta \big)
\big(A{( \Phi^{-1}(z,\omega),\omega)} {\big( \nabla_{\!\! z} + 2i\pi \theta \big)} [\cdot] \big)
\\ 
&+V {\big(\Phi^{-1}\left(z,\omega\right),\omega\big)} [\cdot].
\end{aligned}
\end{equation}
Here, $\Phi: \R^n \times \Omega \to \R^n$ is a stochastic diffeomorphisms,
(called stochastic deformations). 
The $n \times n$ matrix-value stationary function $A = (A_{k \ell})$ and the real stationary potential $V$ are 
measurable and
bounded functions, i.e.
$A_{k \ell}, V \in L^\infty(\mathbb{R}^n \times \Omega)$.
Moreover, $A$ is symmetric and 
uniformly positive defined, that is, there exists $a_0> 0$, such that, for a.a. 
$(y, \omega) \in \mathbb{R}^n \times \Omega$, and each $\xi \in \R^n$
\begin{equation}
\label{ASSUM2}
	\sum_{k,\ell=1}^n A_{k\ell}(y,\omega)\,  \xi_k \, \xi_\ell \geqslant a_0 {\vert \xi \vert}^2.
\end{equation}
The stationarity property of random 
functions will be precisely defined in Section \ref{notation}, also the definition of 
stochastic deformations which were introduced by 
X. Blanc, C. Le Bris, P.-L. Lions (see \cite{BlancLeBrisLions1,BlancLeBrisLions2}).

\medskip
One observes that, the eigenvalue problem \eqref{92347828454trfhfd4rfghjls} 
is called Bloch or shifted spectral cell equation, see Section 2.4 in \cite{VCWNJS}.    
Moreover, each $\theta \in \R^n$ is called a Bloch frequency, $\lambda(\theta)$ is called a Bloch energy and the corresponded 
$\Psi(\theta)$ is called a Bloch wave.

\medskip
Now, let us consider the following spaces
\begin{equation}
	\mathcal{L}_\Phi := {\big\{F(z,\omega)= f( \Phi^{-1} (z, \omega), \omega); f \in L^2_{\rm loc}(\mathbb{R}^n; L^2(\Omega)) \;\; \text{stationary} \big\}}
\end{equation}
and
\begin{equation}
\label{SPACEHPHI}
		\mathcal{H}_\Phi := {\big\{F(z,\omega)= f( \Phi^{-1} (z, \omega), \omega); \; f\in H^1_{\rm loc}(\mathbb{R}^n; L^2(\Omega)) \;\; \text{stationary} \big\}}
\end{equation}
which are Hilbert spaces, endowed respectively with the inner products  
$$
\begin{aligned}
  {\langle F, G \rangle}_{\mathcal{L}_\Phi}&:= \int_\Omega \int_{\Phi([0,1)^n,\omega)} \!\! F(z, \omega) \, \overline{ G(z, \omega) } \, dz \, d\mathbb{P}(\omega),
\\[5pt]
{\langle F, G \rangle}_{\mathcal{H}_\Phi}&:= \int_\Omega \int_{\Phi([0,1)^n,\omega)} \!\! F(z, \omega) \, \overline{ G(z, \omega) } \, dz \, d\mathbb{P}(\omega)
\\
&\; \quad +\int_\Omega \int_{\Phi([0,1)^n,\omega)} \!\! \nabla_{\!\! z} F(z, \omega) \cdot \overline{ \nabla_{\!\! z} G(z, \omega) } \, dz \, d\mathbb{P}(\omega). 
\end{aligned}
$$
\begin{remark}
\label{REMFPHI}
Under the above notations, 
when $\Phi= Id$ we denote $\mathcal{L}_\Phi$ and $\mathcal{H}_\Phi$ by  $\mathcal{L}$ and $\mathcal{H}$ respectively.
Moroever, a function $F \in \clg{H}_\Phi$ if, and only if, $F \circ \Phi \in \clg{H}$.
Analogously, $F \in \clg{L}_\Phi$ if, and only if, $F \circ \Phi \in \clg{L}$. 
\end{remark}

\medskip
Associated to the eigenvalue-eigenfunction problem \eqref{92347828454trfhfd4rfghjls}
is a variational formulation. 
To this end, we consider for $F,G \in \mathcal{H}_\Phi$,  the following functional 
\begin{equation}
\label{FORMVARIAC}
\begin{aligned}
&\langle L^\Phi(\theta)[F], G\rangle
\\
&= \int_\Omega \int_{\Phi([0,1)^n,\omega)} \!\!\!\!\!\!\!\!\! A( \Phi^{-1}(z, \omega), \omega) (\nabla_{\!\! z} + 2i\pi\theta)  F(z,\omega) \cdot 
 \overline{{( \nabla_{\!\! z} + 2i\pi\theta)} G(z,\omega)} \, dz \, d\mathbb{P}(\omega) 
 \\[5pt]
 &+ \int_\Omega \int_{\Phi([0,1)^n,\omega)} V{( \Phi^{-1}(z, \omega), \omega)} \ F(z,\omega) \, 
 \overline{G(z,\omega)} \, dz \, d\mathbb{P}(\omega). 
\end{aligned}
\end{equation}
Then, the variational problem is to minimize $\clg{E}(\Psi):= \langle L^\Phi(\theta)[\Psi], \Psi\rangle$,
subject to the constraint $\|\Psi\|^2_{\mathcal{L}_\Phi} \equiv  {\langle \Psi, \Psi \rangle}_{\mathcal{L}_\Phi}= 1$.   
The minimizing function $\Psi_0$, when it exists, satisfies \eqref{92347828454trfhfd4rfghjls} with $\lambda= \lambda_0$. 

\medskip
A standard routine 
to solve this variational problem relies on a
compact embedding.  
Albeit, from mathematical-physical reasons, see \cite{VCWNJS}, 
the parameter $\omega \in \Omega$ in \eqref{92347828454trfhfd4rfghjls},
equivalently \eqref{FORMVARIAC}, can not be fixed, that is to say,
$\lambda_0$ can not depend on $\omega$. Therefore, 
we established an equivalence between $\mathcal{H}_\Phi$ and
the Sobolev space on groups, see Section \ref{kjh876}, 
and then consider a related Rellich-Kondrachov's Theorem. 
Indeed, we establish a compactness argument in Section \ref{927394r6fy7euh73f}, 
which enable us to show that the space $\mathcal{H}_\Phi$ is compactly embedded in $\mathcal{L}_\Phi$,
in order to solve the associated variational problem. 

\medskip
Finally, we stress that, the subject of Sobolev spaces on Abelian locally compact groups,
to the best of our knowledge,
was introduced by P. G\'orka, E. G. Reyes \cite{GorkaReyes}. 

\subsection{Notation and background} 
\label{notation}
We denote by $\mathbb{G}$ the group $\mathbb{Z}^n$ (or $\mathbb{R}^n$), with $n \in \mathbb{N}$.
The set $[0,1)^n$
denotes the unit cube, which is also called the unitary cell and will be used  
as the reference period for periodic functions.
The symbol $\left\lfloor x \right\rfloor$ denotes the 
unique number in $\mathbb{Z}^n$, such that $x - \left\lfloor x \right\rfloor \in [0,1)^n$.

\medskip
Let $U \subset \R^{n}$ be an open set, $p \geqslant  1$, and $s \in \mathbb{R}$.
We denote by 
$L^p(U)$ the set of (real or complex) $p-$summable functions
with respect to the Lebesgue measure. Given a Lebesgue measurable set
$E \subset \R^n$, 
$|E|$ denotes its $n-$dimensional Lebesgue measure.
Moreover, we will use the standard notations for the 
Sobolev spaces $W^{s,p}(U)$ and $H^{s}(U)\equiv W^{s,2}(U)$. 

\medskip
Now, let $(\Omega, \mathcal{F}, \mathbb{P})$ be a probability space. For each random variable  
$f$ in $L^1(\Omega; \PM)$, ($L^1(\Omega)$ for short), 
we define its expectation value by
$$
    \mathbb{E}[f]:= \int_\Omega f(\omega) \ d\PM(\omega).
$$
A mapping $\tau: \mathbb{G} \times \Omega \to \Omega$ is said a $n-$dimensional dynamical 
system, when
\begin{enumerate}
\item[(i)](Group Property) $\tau(0,\cdot)=id_{\Omega}$ and $\tau(x+y,\omega)=\tau(x,\tau(y,\omega))$ for each $x,y \in  \mathbb{G}$ 
and $\omega\in\Omega$.
\item[(ii)](Invariance) The mappings $\tau(x,\cdot):\Omega\to \Omega$ are $\PM$-measure preserving, that is, for all $x \in  \mathbb{G}$ and 
every $E\in \mathcal{F}$, it follows that
$$
\tau(x,E)\in \mathcal{F},\qquad \PM(\tau(x,E))=\PM(E).
$$
\end{enumerate}
We shall use $\tau(k)\omega$ to denote $\tau(k,\omega)$, and it is usual to say that 
$\tau(k)$ is a discrete (continuous) dynamical system if $k \in \Z^n$ ($k \in \R^n$), although we only stress 
this when it is not obvious from the context. 

\medskip
A measurable function $f$ on $\Omega$ is called $\tau$-invariant, when for each $k \in \mathbb{G}$ 
$$
    f(\tau(k) \omega)= f(\omega) \quad \text{for almost all $\omega \in \Omega$}. 
$$
Then, a measurable set $E \in \mathcal{F}$ is $\tau$-invariant, if its characteristic function $\chi_E$ is $\tau$-invariant. 
It is a straightforward to show that, a $\tau$-invariant set $E$ can be equivalently defined by 
$$
    \tau(k) E= E \quad \text{for each $k \in \mathbb{G}$}.
$$
We say that the dynamical system $\tau$ is ergodic, when
all $\tau$-invariant sets $E$ have measure $\PM(E)$ of either zero or one. 
Equivalently, a dynamical system is ergodic if 
each $\tau$- invariant function is constant almost everywhere, that is 
$$
    \Big( f(\tau(k) \omega)= f(\omega) \quad \text{for each $k \in \mathbb{G}$ and a.e. $\omega \in \Omega$} \Big) 
    \Rightarrow \text{ $f(\cdot)= const.$ a.e.}.  
$$

\medskip
Now, let $(\Gamma, \mathcal{G}, \mathbb{Q})$ be a given probability space. We say that a
measurable function $g: \R^n \times\Gamma \to \R$ is stationary, if for any finite set 
consisting of points $x_1,\ldots,x_j\in \R^n$, and any $k \in \mathbb{G}$, the distribution of the random vector 
$$
   \Big(g(x_1+k,\cdot),\cdots,g(x_j+k,\cdot)\Big)
$$
is independent of $k$. Moreover, subjecting the stationary function $g$ to some natural conditions
it can be showed that, there exists other probability space $(\Omega, \mathcal{F}, \PM)$,  a $n-$dimensional dynamical system 
$\tau: \mathbb{G} \times \Omega \to \Omega$ and a measurable function $f: \R^n \times \Omega \to \R$ satisfying 
\begin{itemize}
\item[(i)] For all $x \in \R^n$, $k \in \mathbb{G}$ and $\PM-$almost every $\omega \in \Omega$ 
\begin{equation}
\label{Stationary}
   f(x+k, \omega)= f(x, \tau(k) \omega).
\end{equation} 

\item[(ii)] For any $x \in \R^n$ the random variables $g(x,\cdot)$, $f(x,\cdot)$ have the same 
law. 
\end{itemize} 

\begin{remark} The set of stationary functions forms an algebra, and
also is stable by limit process. 
Therefore, the product of two
stationaries functions is a stationary one, and the derivative of a 
stationary function is stationary. 
Moreover, the stationarity concept is the most general extension of the 
notions of periodicity and almost periodicity for a function to have some "self-averaging" behaviour. 
\end{remark}

\medskip
To follow, we present the definition of the stochastic deformation as presented in \cite{AndradeNevesSilva}.
\begin{definition}
\label{GradPhiStationary}
A mapping $\Phi: \R^n \times \Omega \to \R^n, (y,\omega) \mapsto z= \Phi(y,\omega)$, is called a stochastic deformation (for short $\Phi_\omega$), when satisfies:
\begin{itemize}
\item[i)] For $\mathbb{P}-$almost every $\omega \in \Omega$, $\Phi(\cdot,\omega)$ is a bi--Lipschitz diffeomorphism.

\item[ii)] There exists $\nu> 0$, such that
$$
\underset{\omega \in \Omega, \, y \in \R^n}{\rm ess \, inf} 
\big({\rm det} \big(\nabla \Phi(y,\omega)\big)\big) \geq \nu.
$$
\item[iii)] There exists a $M> 0$, such that
$$
 \underset{\omega \in \Omega, \, y \in \R^n}{\rm ess \, sup}\big(|\nabla \Phi(y,\omega)|\big) \leq M< \infty.
$$
\item[iv)]
The gradient of $\Phi$, i.e. $\nabla\Phi(y,\omega)$, is stationary in the sense~\eqref{Stationary}.
\end{itemize}
\end{definition}

\medskip
Connected with the notion of stationarity,
we consider now the concept of mean value. 
A function $f \in L^1_{\loc}(\R^n)$ is said to possess a mean value if the 
sequence $\{f(\cdot/\ve){\}}_{\ve>0}$ converges in the duality with $L^{\infty}$ and compactly supported 
functions to a constant $M(f)$. This convergence is equivalent to
\begin{equation}
\label{MeanValue}
\lim_{t\to\infty}\frac1{t^n|A|}\int_{A_t}f(x)\,dx=M(f),
\end{equation}
where $A_t:=\{x\in\R^n\,:\, t^{-1}x\in A\}$, for $t>0$ and any $A \subset \R^n$, with $|A| \ne0$.

\begin{remark}
\label{REMERG}
Unless otherwise stated, we assume that the dynamical system $\tau: \mathbb{G} \times \Omega\to\Omega$ is ergodic 
and we will also use the notation 
$$
   \Medint_{\R^n} f(x) \ dx \quad \text{for $M(f)$}.
$$
\end{remark}
Then, we state the result due to Birkhoff, 
see \cite{Krengel}. 
\begin{theorem}[Birkhoff Ergodic Theorem]\label{Birkhoff}
Let $f \in L^1_\loc(\R^n; L^1(\Omega))$ be a stationary random variable. 
Then, for almost every $\widetilde{\omega} \in \Omega$ the function 
$f(\cdot,\widetilde{\omega})$ possesses a mean value in the sense of~\eqref{MeanValue}. Moreover, the mean value 
$M\left(f(\cdot,\widetilde{\omega})\right)$ as a function of $\widetilde{\omega} \in\Omega$ satisfies
for almost every $\widetilde{\omega} \in \Omega$: 

\smallskip
i) Discrete case (i.e. $\tau: \Z^n \times \Omega \to \Omega$);
$$
    \Medint_{\R^n} f(x,\widetilde{\omega}) \ dx= 
   \mathbb{E} \left[\int_{[0,1)^n}  f(y,\cdot)\, dy\right].
$$

ii) Continuous case (i.e. $\tau: \R^n \times \Omega \to \Omega$);
$$
    \Medint_{\R^n} f(x,\widetilde{\omega}) \ dx=  \mathbb{E}\left[ f(0,\cdot) \right].
$$
\end{theorem}

\medskip
The Birkhoff Ergodic Theorem holds if a stationary function is composed 
with a stochastic deformation: 
\begin{lemma}\label{phi2}
Let $\Phi$ be a stochastic deformation and $f \in L^{\infty}_\loc(\R^n; L^1(\Omega))$ be a stationary random variable in the 
sense~\eqref{Stationary}. Then for almost $\widetilde{\omega} \in \Omega$ the function 
$f\left(\Phi^{-1}(\cdot,\widetilde{\omega}),\widetilde{\omega}\right)$ possesses a mean value 
in the sense of~\eqref{MeanValue} and satisfies: 

\smallskip
i) Discrete case;
$$
\text{$\Medint_{\R^n}f\left(\Phi^{-1}(z,\widetilde{\omega}),\widetilde{\omega}\right)\,dz
= \frac{\mathbb{E}\left[\int_{\Phi([0,1)^n, \cdot)} f {\left( \Phi^{-1}\left( z, \cdot \right), \cdot \right)} \, dz \right]}
{\det\left(\mathbb{E}\left[\int_{[0,1)^n} \nabla_{\!\! y} \Phi(y,\cdot) \, dy \right]\right)}$
\quad for a.a. $\widetilde{\omega} \in \Omega$}.
$$

ii) Continuous case; 
$$
\text{$\Medint_{\R^n}f\left(\Phi^{-1}(z,\widetilde{\omega}),\widetilde{\omega}\right)\,dz
= \frac{\mathbb{E}\left[f(0,\cdot)\det\left(\nabla\Phi(0,\cdot)\right)\right]}
{\det\left(\mathbb{E}\left[\nabla \Phi(0,\cdot)\right]\right)}$
\qquad for a.a. $\widetilde{\omega} \in \Omega$}.
$$
\end{lemma}

\begin{proof}
See Blanc, Le Bris, Lions \cite{BlancLeBrisLions1}, also
Andrade, Neves, Silva \cite{AndradeNevesSilva}.
\end{proof}

The next theorem presents important properties of
stationary functions. 
\begin{theorem}
\label{987987789879879879}
For $p> 1$, let $u,v \in L^1_{\rm loc}(\mathbb{R}^n; L^p(\Omega))$
be stationary functions. Then, for any $i \in \{1,\ldots,n \}$ fixed, the following sentences are equivalent: 
\begin{equation}
\label{837648726963874}
(A) \quad  \int_{[0,1)^n} \int_\Omega u(y,\omega) \frac{\partial {\zeta}}{\partial y_i} (y, \omega) \, d\mathbb{P}(\omega) \, dy 
      = - \int_{[0,1)^n} \int_\Omega v(y,\omega) \, {\zeta}(y,\omega) \, d\mathbb{P} \, dy, \hspace{20pt}
\end{equation}
for each stationary function $\zeta \in C^1( \mathbb{R}^n; L^q(\Omega))$,
with $1/p + 1/q = 1$. 
\begin{equation}
\label{987978978956743}
(B) \quad  \int_{\mathbb{R}^n} u(y,\omega) \frac{\partial {\varphi}}{\partial y_i} (y) \, dy = - \int_{\mathbb{R}^n} v(y,\omega) \, {\varphi}(y) \, dy,
\hspace{87pt}
\end{equation}
for any $\varphi \in C^1_{\rm c}(\mathbb{R}^n)$, and almost sure $\omega \in \Omega$.
\end{theorem}

\begin{proof}
See Blanc, Le Bris, Lions \cite{BlancLeBrisLions1}. 
\end{proof} 

\section{Sobolev spaces on groups}
\label{9634783yuhdj6ty}

To begin, we sum up some definitions and properties of topological groups, 
which will be used along this section. Most of the material could be found in 
E. Hewitt, A. Ross \cite{HewittRoss} and G. B. Folland \cite{Folland2}
(with more details). 

\medskip
A nonempty set $G$ endowed with an application, $\ast : G \! \times \! G \to G$,
is called a group, when for each $x, y, z \in G$:
\begin{itemize}
\item[1.] ${ (x\ast y) \ast z = x \ast (y \ast z) }$;
\item[2.] There exists ${e \in G }$, such that ${ x \ast e = e \ast x = e }$;
\item[3.] For all ${ y \in G }$, there exists ${y^{-1}\in G }$, such that ${ y \ast y^{-1} = y^{-1} \ast y = e }$.
	\end{itemize}
Moreover, if $x \ast y = y \ast x$, then $G$ is called an Abelian group. 
From now on, we write for simplicity $x \, z$ instead of $x \ast z$. 
A topological group is a group $G$  together with a topology, such that,
both the group's binary operation $(x,y) \mapsto x \, y$,
and the function mapping group elements to their respective inverses 
$x \mapsto x^{-1}$
are continuous functions with respect to the topology.
Unless the contrary is explicit stated, any group mentioned here is 
a locally compact Abelian (LCA for short) group, and 
we may assume without loss of generality that, 
the associated topology is Hausdorff 
(see G. B. Folland \cite{Folland2}, Corollary 2.3).

\medskip
A complex value function
$\xi : G \to \mathbb{S}^1$ is called a character of $G$, when 
$$
 \xi(x \, y) = \xi(x) \xi(y),
\quad  \quad  \text{(for each $x, y \in G$)}.
 $$
 We recall that, the set of characters of $G$ is an Abelian group
 with the usual product of functions, identity element $e= 1$, and
 inverse element $\xi^{-1} = \overline{\xi}$.
 The characters' group of the topological group $G$, called
the dual group of $G$ and denoted by $G^\wedge$,  
is the set of all continuous characters, that is to say 
$$
 G^\wedge:= \{ \xi : G \to \mathbb{S}^1 \; ; \; \text{$\xi$ is a continuous homomorphism}\}.
$$
Moreover, we may endow $G^\wedge$ with a topology with respect to which,
$G^\wedge$ itself is a LCA group.  

\medskip
We denote by $\mu$, $\nu$ the unique (up to a positive multiplicative constant) Haar mesures in $G$ and $G^\wedge$ respectively. 
The $L^p$ spaces over $G$ and its dual are defined as usual, with their respective mesures. 
Let us recall two important properties when $G$ is compact:
\begin{equation}
\label{CARACGCOMP}
\begin{aligned}
&i) \quad \text{If $\mu(G)= 1$, then $G^\wedge$ is an orthonormal set in $L^2(G;\mu)$}.
\\[5pt]
&ii) \quad \text{The dual group $G^\wedge$ is discrete, and $\nu$ is the countermeasure}. 
\end{aligned}
\end{equation}

\medskip
One remarks that, the study of Sobolev 
spaces on LCA groups uses essentially the concept of Fourier Transform, then we have the following   
\begin{definition}
Given a complex value function $f \in L^1(G;\mu)$, the function $\widehat{f}: G^\wedge \to \mathbb{C}$, defined by
\begin{equation}
	\widehat{f}(\xi):= \int_G f(x) \, \overline{\xi(x)} \, d\mu(x)
\end{equation}	
is called the Fourier transform of $f$ on $G$.
\end{definition}
Usually, the Fourier Transform of $f$ is denoted by $\clg{F}f$ to emphasize that it is an operator, 
but we prefer to adopt the usual notation $\widehat{f}$. 
Moreover, we recall that the Fourier transform is an homomorphism from $L^1(G;\mu)$ to $C_0(G^\wedge)$ 
(or $C(G^\wedge)$ when $G$ is compact), see Proposition 4.13 in \cite{Folland2}. Also we address the reader to 
 \cite{Folland2}, Chapter 4, for the Plancherel Theorem
and the Inverse Fourier Transform. 

\medskip
Before we establish the definition of (energy) Sobolev spaces on LCA groups, let us
consider the following set
$$
\begin{aligned}
  {\rm P}= \{p: G^\wedge \times & G^\wedge \to [0,\infty) / 
  \\
  \; & \text{$p$ is a continuous invariant pseudo-metric on $G^\wedge$} \}.
\end{aligned} 
$$
The Birkhoff-Kakutani Theorem (see \cite{HewittRoss} p.68) 
ensures that, the set P is not empty. 
Any pseudo-metric $p \in {\rm P}$ is well defined for each $(x,y) \in G^\wedge \times G^\wedge$, hence we may define
\begin{equation}
\label{Gamma}
\gamma(x):= p(x,e) \equiv p(x,1).
\end{equation} 
Moreover, one observes that $\gamma(1)= 0$.
Then, we have the following 
\begin{definition}[Energy Sobolev Spaces on LCA Groups]
\label{SOBOLEVESPACES}
Let $s$ be a non-negative real number and $\gamma(x)$ be given by \eqref{Gamma}
for some fixed $p \in {\rm P}$. The energy Sobolev space 
$H^s_\gamma(G)$ is the set of functions $f \in L^2(G;\mu)$, such that
\begin{equation}
     \int_{G^\wedge} (1+\gamma(\xi)^2)^s  \, |\widehat{f}(\xi)|^2 d\nu(\xi)< \infty.
\end{equation}
Moreover, given a function $f \in H^s_\gamma(G)$ its norm is defined as 
\begin{equation}
   \Vert f \Vert_{H^s_\gamma(G)} := \left( \int_{G^\wedge} \left(1+\gamma(\xi)^2 \right)^s \vert \widehat{f}(\xi) \vert^2 d\nu(\xi) \right)^{1/2}.
\end{equation}
\end{definition}
Below, taking specific functions $\gamma$, the usual Sobolev spaces on $\R^d$ and 
other examples are considered. In particular,   
the Plancherel Theorem implies that, $H^0_\gamma(G)=  L^2(G;\mu)$.

\begin{example}
\label{EXAMPLERN}
Let $G= (\R^n, +)$ which is known to be a LCA group, and consider 
its dual group $(\mathbb{R}^n)^\wedge = \{ \xi_y \; ; \; y\in\mathbb{R}^n \}$,
where for each $x \in \R^n$
\begin{equation}
\label{caracterunitario}
  \xi_y(x) = e^{2 \pi i \, y \cdot x},
\end{equation}
hence $|\xi_y(x)|= 1$ and $\xi_0(x)= 1$. 
One remarks that, here we denote (without invocation of vector space structure)
$$
    a \cdot b= a_1 b_1 + a_2 b_2 + \ldots + a_n b_n, \quad \text{(for all $a,b \in G$)}.
$$
For any $x, y \in \R^n$ let us consider 
$$
  p(\xi_x,\xi_y)= 2\pi  \|x - y\|, 
$$
where $\| \cdot \|$ is the Euclidean norm in $\R^n$. Hence $\gamma(\xi_x)= p(\xi_x,1)= 2 \pi \|x\|$.  
Since $(\mathbb{R}^n)^\wedge \cong \mathbb{R}^n$, the Sobolev space $H^s_\gamma(G)$ coincide 
with the usual Sobolev space on $\R^n$. 
\end{example}

\begin{example}
\label{6576dgtftdefd}
Let us recall that, the set $[0,1)^n$ endowed with the binary
operation 
$$
   (x,y) \in [0,1)^n \! \times \! [0,1)^n \;\; \mapsto \;\; x+y - \left\lfloor x+y \right\rfloor \in [0,1)^d
$$ 
is an Abelian group, and the function 
$\Lambda: \mathbb{R}^n \to [0,1)^n$, $\Lambda(x):= x - \left\lfloor x \right\rfloor$
is an homomorphism of groups. Moreover, under the 
induced topology by $\Lambda$, that is to say 
$$
  \{U \subset [0,1)^n \; ; \; \Lambda^{-1}(U) \; \text{is an open set of} \;\, \mathbb{R}^n \}, 
$$
 $[0,1)^n$ is a compact Abelian group, which is called $n-$dimensional Torus and denoted by 
$\mathbb{T}^n$. Its dual group is characterized by the integers $\Z^n$, that is 
$$
\text{
$(\mathbb{T}^n)^\wedge = \{ \xi_m \; ; \; m \in \mathbb{Z}^n \}$, where $\xi_m(x)$ is given by  
\eqref{caracterunitario} for all $x \in \mathbb{R}^n$}. 
$$
For each $m,k \in \Z^n$, we consider 
$$
   p(\xi_m,\xi_k)= 2\pi \sum_{j=1}^n{\vert m_j - k_j \vert},
   \quad \text{and thus $\gamma(\xi_m)= 2 \pi \sum_{j=1}^n{\vert m_j \vert}$}.
$$
Then, the Sobolev space $H^s_\gamma(\mathbb{T}^n)$ coincide 
with the usual Sobolev space on $\mathbb{T}^n$.
		
\smallskip
Now, following the above discussion let us consider the infinite Torus 
$\mathbb{T}^I$, where $I$ is an index set. Since an arbitrary product of compact spaces is compact in the 
product topology (Tychonoff Theorem), $\mathbb{T}^I$ is a compact Abelian group. Here, 
the binary operation on $ \mathbb{T}^I \times  \mathbb{T}^I$ is defined coordinate by coordinate, that is, for each 
$\ell \in I$ 
$$
   g_\ell + h_\ell:= g_\ell + h_\ell - \left\lfloor g_\ell + h_\ell \right\rfloor.
$$
Moreover, the set $\mathbb{Z}^I_{\rm c} := \{ m \in \mathbb{Z}^I; \text{{\rm supp} $m$ is compact} \}$
characterizes the elements of the dual group $(\mathbb{T}^I)^\wedge$.
 Indeed, applying Theorem 23.21 in 
\cite{HewittRoss}, similarly we have  
$$
   (\mathbb{T}^I)^\wedge = {\left\{ \xi_m \; ; \; m\in \mathbb{Z}^I_{\rm c} \right\}},
$$
where $\xi_m(k)$ is given by  
\eqref{caracterunitario} for each $m,k \in \mathbb{Z}_{\rm c}^I$, the pseudo-metric
$$
   p(\xi_m,\xi_k)= 2\pi \sum_{\ell \in I}{\vert m_\ell - k_\ell \vert},
   \quad \text{and $\gamma(\xi_m)= 2 \pi \sum_{\ell \in I}{\vert m_\ell \vert}$}.
$$
Consequently, we have establish the Sobolev spaces $H^s_{\gamma}(\mathbb{T}^I)$.
\end{example}

\section{Groups and Dynamical systems}
\label{kjh876}

In this section, we are interested to come together the discussion 
about dynamical systems considered in Section \ref{notation}
with the theory developed in the last section 
for LCA groups. To this end, we consider 
stationary functions in the continuous sense (continuous dynamical systems). 
Moreover, we recall that all the groups in this paper are 
assumed to be Hausdorff.  

\medskip
To begin, let $G$ be a locally compact group with Haar measure $\mu$,
we know that $\mu(G)< \infty$ if, and only if, $G$ is compact. 
Therefore, we consider from now on that $G$ is a compact Abelian group, 
hence $\mu$ is a finite measure and, up to a normalization, $(G,\mu)$ is a probability space.
%
%
We are going to consider the dynamical systems, $\tau: \R^n \times G \to G$, defined by 
\begin{equation}
\label{TAUFI}
   \tau(x) \omega:= \varphi(x) \, \omega,
\end{equation}
where $\varphi: \R^n \to G$ is a given (continuous) homomorphism. 
Indeed, first $\tau(0) \omega= \omega$ and 
$\tau(x+y, \omega)= \vp(x) \vp(y) \omega= \tau(x,\tau(y)\omega)$. 
Moreover, since $\mu$ is a translation invariant Haar measure, the 
mapping $\tau(x,\cdot): G \to G$ is $\mu-$measure preserving. 
Recall from Remark \ref{REMERG} we have assumed that, 
the dynamical systems we are interested here are
ergodic. Then, it is important to characterize the conditions
for the mapping $\vp$, under which the dynamical system defined by 
\eqref{TAUFI} is ergodic. To this end, first let us consider the following 

\begin{lemma}
\label{DIST}
Let $H$ be a topological group, $F \subset H$ closed, $F \neq H$ and $x \notin F$.
Then, there exists a neighborwood $V$ of the identity $e$, such that
$$
   F V \cap x V= \emptyset. 
$$
\end{lemma}

\begin{proof}
First, we observe that:

i) Since $F \subset H$ is closed and $F \neq H$, there
exists a neighborwood $U$ of the identity $e$,
such that $F \cap x U= \emptyset$. 

ii) There exists a symmetric neighborwood $V$ of the identity $e$,
such that $VV \subset U$. 

Now, suppose that $F V \cap x V \neq \emptyset$. 
Therefore, there exist $v_1, v_2 \in V$ and $k_0 \in F$ such that, $k_0 v_1= x v_2$.   
Consequently, $k_0= x v_2 v_1^{-1}$ and from $(ii)$ it follows that, $k_0 \in x U$. 
Then, we have a contradiction from $(i)$. 
\end{proof}

 \underline {\bf Claim 1:} The dynamical system defined 
by \eqref{TAUFI} is ergodic if, and only if, 
$\vp(\R^n)$ is dense in $G$. 

\smallskip
Proof of Claim 1: 1. Let us show first the necessity. Therefore, we suppose that 
$\vp(\R^n)$ is not dense in $G$, that is $K:= \overline{\vp(\R^n)} \neq G$. 
Then, applying Lemma \ref{DIST}
there exists a neighborhood $V$ of $e$, such that $K V \cap x V= \emptyset$,
for some $x \notin K$. Recall that the Haar measure on open sets are positive, 
moreover
$$
   K V= \bigcup_{k \in K} k V,
$$
which is an open set, thus we have 
$$
   0< \mu(K V) + \mu(x V) \leq 1. 
$$
Consequently, it follows that $0< \mu(\vp(\R^n) V)< 1$. For convenience, le us denote 
$E= \vp(\R^n) V$, hence $\tau(x) E= E$ for each $x \in \R^n$. 
Then, the dynamical system $\tau$ is not ergodic, since $E \subset G$ is a $\tau$-invariant set 
with $0< \mu(E)< 1$. 

\medskip
2. It remains to show the sufficiency. 
Let $E \subset G$ be a $\mu-$measurable $\tau$-invariant set,
hence $\omega E= E$ for each $\omega \in \vp(\R^n)$. Assume 
by contradiction that, $0< \mu(E)< 1$, thus $\mu(G \setminus E)> 0$.
Denote by $\mathcal{B}_G$ the Borel $\sigma-$algebra on $G$, and define, 
$\lambda:= \mu_{\lfloor E}$, that is $\lambda(A)= \mu(E \cap A)$ for all 
$A \in \mathcal{B}_G$. Recall that $G$ is not necessarily metric, therefore, it is not
clear if each Borel set is $\mu-$measurable. Then, it follows that: 

$(i)$ For any $A \in \mathcal{B}_G$ fixed, the mapping 
$\omega \in G \mapsto \lambda(\omega A)$ is continuous. 
Indeed, for $\omega \in G$ and $A \in \mathcal{B}_G$, we have
$$
\begin{aligned}
 \lambda(\omega A)&= \int_G 1_E(\varpi) 1_{\omega A}(\varpi) d\mu(\varpi)
 \\[5pt]
 &= \int_G 1_E(\varpi) 1_{A}(\omega^{-1} \varpi) d\mu(\varpi) 
= \int_G 1_E(\omega \varpi) 1_{A}(\varpi) d\mu(\varpi).
\end{aligned}
$$
Therefore, for $\omega, \omega_0 \in G$
$$
\begin{aligned}
|\lambda(\omega A) - \lambda(\omega_0 A)|&= \big| \int_G \big(1_E(\omega \varpi) - 1_E(\omega_0 \varpi)\big) 1_A(\varpi)  d\mu(\varpi) \big|
\\
 &\leq \big(\mu(A)\big)^{1/2}
 \big( \int_G  |1_E(\omega \varpi) -  1_E(\omega_0 \varpi)|^2 d\mu(\varpi) \big)^{1/2}
  \tobo{\omega \to \omega_0} 0. 
\end{aligned}
$$

$(ii)$ $\lambda$ is invariant, i.e. for all $\omega \in G$, and $A \in \mathcal{B}_G$, $\lambda(\omega A)= \lambda(A)$. 
Indeed, for each $\omega \in \vp(\R^d)$, and $A \in \mathcal{B}_G$, we have 
$$
   (\omega A) \cap E= (\omega A) \cap (\omega E)= \omega (A \cap E).  
$$
Thus since $\mu$ is invariant, $\mu_{\lfloor E}(\omega A)= \mu_{\lfloor E}(A)$. Consequently,
due to item $(i)$ and $\overline{\vp(\R^d)}= G$, it follows that $\lambda$ is invariant. 

From item $(ii)$ the Radon measure $\lambda$ is a Haar measure on $G$. By the uniqueness 
of the Haar measure on $G$, there exists $\alpha> 0$, such that for all $A \in \mathcal{B}_G$,
$\alpha \lambda(A)= \mu(A)$. In particular, $\alpha \lambda(G \setminus E)= \mu(G \setminus E)$.
But $\lambda(G \setminus E)= 0$ by definition and $\mu(G \setminus E)> 0$, which is a contradiction
and hence $\tau$ is ergodic. 

\begin{remark}
1. One remarks that, in order to show that $\tau$ given by \eqref{TAUFI} is ergodic, it was not used
  that $\vp$ is continuous, nor that $G$ is metric. Compare with the statement in \cite{JikovKozlovOleinik} 
  p.225 (after Theorem 7.2). 

2. From now on, we assume that $\vp(\R^n)$ is dense in $G$. 
\end{remark}

\medskip
Now, for the dynamical system established before, the main issue is to show how the Sobolev space 
$H^1_{\gamma}(G)$ is related with the space $\mathcal{H}_\Phi$ given by \eqref{SPACEHPHI} for $\Phi= Id$, 
that is 
$$
    \mathcal{H}= {\big\{f(y, \omega); \; f \in H^1_{\rm loc}(\mathbb{R}^n; L^2(G)) \;\; \text{stationary} \big\}},
$$
which is a Hilbert space endowed with the following inner product  
$$
{\langle f,g \rangle}_{\mathcal{H}}= \int_G  f(0, \omega) \, \overline{g(0, \omega) } \, d\mu(\omega)
+ \int_G \nabla_{\!\! y} f(0, \omega) \cdot \overline{ \nabla_{\!\! y} g(0, \omega) }\, d\mu(\omega).
$$
Let $\chi$ be a character on $G$, i.e. $\chi \in G^\wedge$. Since $\vp: \R^n \to G$ is a continuous homomorphism, the 
function $(\chi \circ \vp): \R^n \to \C$
is a continuous character in $\R^n$. More precisely, given any fixed $\chi \in G^\wedge$ we may find 
$y \in \R^n$, $(y \equiv y(\chi))$, such that, for each $x \in \R^n$
$$
  \big(\chi \circ \vp \big)(x) =:\xi_{y(\chi)}(x)= e^{2\pi i \, y(\chi) \cdot x}.
$$
Following Example \ref{EXAMPLERN} we define the pseudo-metric 
$p_\vp: G^\wedge \times G^\wedge \to [0,\infty)$ by  
\begin{equation}
\label{PSEDO}
   p_\vp(\chi_1, \chi_2):= p(\xi_{y_1(\chi_1)}, \xi_{y_2(\chi_2)})= 2 \pi \|y_1(\chi_1) - y_2(\chi_2)\|. 
\end{equation}
Then, we have 
$$
   \gamma(\chi)= p_\vp(\chi,1)=  2 \pi \|y(\chi)\|. 
$$

\medskip
Let us observe that, we have used in the above construction of $\gamma$ the continuity of the homomorphism $\vp: \R^n \to G$,
that is to say, it was essential the continuity of $\vp$. In fact, the function $\gamma$ was given by the pseudo-metric $p_\vp$, which is 
not necessarily a metric. Although, we have the following 

\medskip
 \underline {\bf Claim 2:} The pseudo-metric $p_\vp: G^\wedge \times G^\wedge \to [0,\infty)$ given by \eqref{PSEDO} 
is a metric if, and only if, $\vp(\R^n)$ is dense in $G$. 

\smallskip
Proof of Claim 2: 1. First, let us assume that $\overline{\vp(\R^n)} \neq G$, and then show that $p_\vp$ is not a metric. 
Therefore, we have the necessity proved. From Corollary 24.12 in \cite{HewittRoss}, since 
$\overline{\vp(\R^n)}$ is a closer proper subgroup of $G$, hence there exists $\xi \in G^\wedge \setminus \{1\}$,
such that $\xi(\overline{\vp(\R^n)})= \{1\}$. Hence  there exists $\xi \in G^\wedge \setminus \{1\}$, 
such that, $\xi(\vp(x))= 1$, for each $x \in \R^n$, i.e. $y(\xi)= 0$. Therefore, we have 
$p_\vp(\xi, 1)= 0$, 
which implies that $p_\vp$ is not a metric. 

\medskip
2. Now, let us assume that $\overline{\vp(\R^n)}= G$, and it is enough to show that
if $p_\vp(\xi, 1)= 0$, then $\xi= 1$. Indeed, if $0= p_\vp(\xi,1)= 2 \pi \|y(\xi)\|$, then $y(\xi)= 0$. 
Therefore, $\xi(\vp(x))= 1$ for each $x \in \R^d$, since $\xi$ is continuous and  $\overline{\vp(\R^n)}= G$,
it follows that, for each $\omega \in G$, $\xi(\omega)= 1$, which finishes the proof of the claim. 

\begin{remark}
Since we have already assumed that $\vp(\R^n)$ is dense in $G$, it follows that 
$p_\vp$ is indeed a metric, which does not imply necessarily that $G$, itself, is metric.   
\end{remark}		

Under the assumptions considered above, we have the following 
\begin{lemma} If $f \in \mathcal{H}$, then for $j \in \{1,\ldots,d\}$ and all $\xi \in G^\wedge$
\begin{equation}
\label{DERIVGROUPFOURIER}
 \widehat{\partial_j f(0,\xi)}= 2 \pi i \; y_j(\xi) \widehat{f(0,\xi)}.
\end{equation}
\end{lemma}

\begin{proof}
First, for each $x \in \R^d$ and $\omega \in G$, define 
$$
\begin{aligned}
   \xi_\tau(x,\omega)&:= \xi(\tau(x,\omega))= \xi(\vp(x) \omega)= \xi(\vp(x)) \; \xi(\omega)
   \\[5pt]
   &= e^{2 \pi i x \cdot y(\xi)} \; \xi(\omega). 
\end{aligned}
$$   
Therefore $\xi_\tau \in C^\infty(\R^d; L^2(G))$, and we have for $j \in \{1,\ldots,d\}$
\begin{equation}
\label{AUXIL}
\partial_j  \xi_\tau(0,\omega)= 2 \pi i \; y_j(\xi) \;  \xi(\omega). 
\end{equation}
Finally, applying Theorem \ref{987987789879879879} we obtain
$$
\begin{aligned}
  \int_G \partial_j f(0,\omega) \; \overline{\xi_\tau}(0,\omega) d\mu(\omega)&= -  \int_G f(0,\omega) \; \partial_j \overline{\xi_\tau}(0,\omega) d\mu(\omega)
  \\[5pt]
 &=  2 \pi i \; y_j(\xi) \int_G f(0,\omega) \; \overline{\xi}(\omega) d\mu(\omega),
\end{aligned} 
$$
where we have used \eqref{AUXIL}. From the above equation and the definition of the 
Fourier transform on groups we obtain \eqref{DERIVGROUPFOURIER}, and the lemma is proved. 
\end{proof}

\medskip
Now we are able to state the equivalence between the spaces $\mathcal{H}$ and $H^1_\gamma(G)$,
which is to say, we have the following 
\begin{theorem}
\label{THMEQNOM}
A function $f \in \mathcal{H}$ if, and only if, $f(0,\cdot) \in H^1_\gamma(G)$,
and 
$$
    \Vert f \Vert_\mathcal{H} = \Vert f(0,\cdot) \Vert_{H_\gamma^1(G)}.
$$
\end{theorem}

\begin{proof}
1. Let us first show that, if $f \in \mathcal{H}$ then $f \in H^1_\gamma(G)$. 
To follow we observe that 
$$
\begin{aligned}
  \int_{G^\wedge} (1 + \gamma(\xi)^2) |\widehat{f(0,\xi)}|^2 \; d\nu(\xi)&= 
   \int_{G^\wedge} |\widehat{f(0,\xi)}|^2 \; d\nu(\xi)
   \\[5pt]
   &+  \int_{G^\wedge} | 2 \pi i \; y(\xi) \widehat{f(0,\xi)}|^2 \; d\nu(\xi)
   \\[5pt]
   &= \int_{G^\wedge} |\widehat{f(0,\xi)}|^2 \; d\nu(\xi)
   +  \int_{G^\wedge} |\widehat{\nabla_{\!\!y} f(0,\xi)}|^2 \; d\nu(\xi),
\end{aligned}
$$
where we have used \eqref{DERIVGROUPFOURIER}. Therefore, applying 
Plancherel theorem 
$$
  \int_{G^\wedge}\! (1 + \gamma(\xi)^2) |\widehat{f(0,\xi)}|^2 \; d\nu(\xi)= \!\!
   \int_{G}\!  |{f(0,\omega)}|^2 \; d\mu(\omega)
   + \! \int_{G} |\nabla_{\!\!y} {f(0,\omega)}|^2 \; d\mu(\omega)\!<  \! \infty,
$$
and thus $f(0,\cdot) \in H^1_\gamma(G)$. 

\medskip
2. Now, let $f(x,\omega)$ be a stationary function, such that $f(0,\cdot) \in H^1_\gamma(G)$, then we show that 
$f \in \mathcal{H}$. Given a stationary function $\zeta \in C^1(\R^d; L^2(G))$, applying the Palncherel theorem and polarization identity
$$
   \int_G \partial_j \zeta(0,\omega) \;  \overline{f(0,\omega)} d\mu(\omega)
  = \int_{G^\wedge} \widehat{\partial_j \zeta(0,\xi)} \; \overline{\widehat{f(0,\xi)}} d\nu(\xi)
$$
for $j \in \{1,\ldots,d\}$. Due to \eqref{DERIVGROUPFOURIER}, we may write
\begin{equation}
\label{HH1}
\begin{aligned}
   \int_G \partial_j \zeta(0,\omega) \;  \overline{f(0,\omega)} d\mu(\omega)
  &= \int_{G^\wedge} 2 \pi
  i \; y_j(\xi)\widehat{\zeta(0,\xi)} \; \overline{\widehat{f(0,\xi)}} d\nu(\xi)
\\[5pt]
&= - \int_{G^\wedge} \widehat{\zeta(0,\xi)} \; \overline{2 \pi i \; y_j(\xi) \widehat{f(0,\xi)}} d\nu(\xi).
\end{aligned}
\end{equation}
For $j \in \{1,\ldots,d\}$ we define, $g_j(\omega):= \big(2 \pi i \; y_j(\xi)  \widehat{f(0,\xi)}\big)^\vee$,
then $g_j \in L^2(G)$. Indeed, we have 
$$
   \int_G |g_j(\omega)|^2 d\mu(\omega)=  \int_{G^\wedge} |\widehat{g_j(\xi)}|^2 d\nu(\xi) 
   \leq \int_{G^\wedge} (1 + \gamma(\xi)^2) |\widehat{f(0,\xi)}|^2 d\nu(\xi)< \infty.
$$
Therefore, we obtain from \eqref{HH1}
$$
   \int_G \partial_j \zeta(0,\omega) \;  \overline{f(0,\omega)} d\mu(\omega)
   = -  \int_G \zeta(0,\omega) \;  \overline{g_j(\omega)} d\mu(\omega)
$$
for any stationary function $\zeta \in C^1(\R^d; L^2(G))$, and $j \in \{1,\ldots,d\}$. 
Then $f \in \mathcal{H}$ due to Theorem \ref{987987789879879879}. 
\end{proof} 
	
\begin{corollary}
Let $f \in L^2_{\loc}(\R^d; L^2(G))$ be a stationary function
and $\Phi$ a stochastic deformation. 
Then, $f \circ \Phi^{-1} \in \clg{H}_\Phi$ if, and only if, $f(0,\cdot) \in H^1_\gamma(G)$,
and there exist constants $C_1, C_2> 0$, such that 
$$
   C_1  \Vert f \circ \Phi^{-1} \Vert_{\mathcal{H}_\Phi}\leq \Vert f(0,\cdot) \Vert_{H_\gamma^1(G)}
   \leq C_2 \Vert f \circ \Phi^{-1} \Vert_{\mathcal{H}_\Phi}.
$$
\end{corollary}
	
\begin{proof}
Follows from Theorem \ref{THMEQNOM} and Remark \ref{REMFPHI}. 
\end{proof}	
	
\section{Rellich--Kondrachov Theorem}
\label{927394r6fy7euh73f}

The aim of this section is to characterize when 
the Sobolev space $H^1_\gamma(G)$ is compactly embedded in $L^2(G)$,
written $H^1_\gamma(G) \subset \subset L^2(G)$, where $G$ is considered a compact Abelian group 
and $\gamma: G^{\wedge} \to [0,\infty)$ is given by \eqref{Gamma}. 
We observe that, $H^1_\gamma(G) \subset \subset L^2(G)$ is exactly the 
Rellich--Kondrachov Theorem on compact Abelian groups, which was established
under some conditions on $\gamma$
in \cite{GorkaReyes}. 
Nevertheless, as a byproduct of the characterization established here, we provide
the proof of this theorem in a more 
precise context. 

\medskip
To start the investigation, let $(G,\mu)$ be a probability space and consider 
the operator
$T: L^2(G^\wedge) \to L^2(G^\wedge)$,
defined by
\begin{equation}
\label{TCOMP}
	[T(f)](\xi) := \frac{f(\xi)}{\sqrt{(1 + \gamma(\xi)^2)}}.
\end{equation}
We remark that, $T$ as defined above is a bounded linear, ($\Vert T \Vert \leqslant 1$),  self-adjoint operator,
which is injective and satisfies for each $f \in L^2(G^\wedge)$
\begin{equation}
\label{76354433}
	\int_{G^\wedge} \left(1 + \gamma(\xi)^2 \right) \, {\vert [T(f)](\xi) \vert}^2 d\nu(\xi) 
	= \int_{G^\wedge} \vert f(\xi) \vert^2 d\nu(\xi).
\end{equation}
Moreover, a function $f \in H^1_\gamma(G)$ if, and only if, $\widehat{f} \in T(L^2(G^\wedge))$, 
that is to say 
\begin{equation}
\label{87648764}
  f \in H^1_\gamma(G) \Leftrightarrow  \widehat{f} \in T(L^2(G^\wedge)). 
\end{equation}
Indeed, if $ f \in H^1_\gamma(G)$ then, we have $f \in L^2(G)$ and 
$$
   \int_{G^\wedge} \left( 1+\gamma(\xi)^2 \right) \vert \widehat{f}(\xi) \vert^2 d\nu(\xi) 
   = \int_{G^\wedge} \vert \sqrt{\left( 1+\gamma(\xi)^2 \right)} \, \widehat{f}(\xi) \vert^2 d\nu(\xi)< \infty.
$$
Therefore, defining $g(\xi):=  \sqrt{\left( 1+\gamma(\xi)^2 \right)} \widehat{f(\xi)}$, hence $g \in L^2(G^\wedge)$ and we have
$\widehat{f} \in T(L^2(G^\wedge))$.

\medskip
Now, if $\widehat{f} \in T(L^2(G^\wedge))$ let us show that, $f \in H^1_\gamma(G)$. First, there exists 
$g \in L^2(G^\wedge)$ such that, $\widehat{f} = T(g)$. 
Thus from equation \eqref{76354433}, we obtain 
$$
    \int_{G^\wedge} (1 + \gamma(\xi)^2) \, |\widehat{f}(\xi)|^2 d\nu(\xi) 
    = \int_{G^\wedge} |g(\xi)|^2 d\nu(\xi)< \infty,
$$
that is, by definition $f \in H^1_\gamma(G)$.

\medskip
Then we have the following Equivalence Theorem:
\begin{theorem}
\label{876876872}
The Sobolev space $H^1_\gamma(G)$ is compactly embedded in $L^2(G)$ 
if, and only if, the operator $T$ defined by \eqref{TCOMP} is compact. 
\end{theorem}

\begin{proof}
1. First, let us assume that $H^1_\gamma(G) \subset \subset L^2(G)$, 
and take a bounded sequence $\{f_m\}$, $f_m \in L^2(G^\wedge)$ 
for each $m \in \N$. Thus $T(f_m) \in L^2(G^\wedge)$, and defining 
$g_m:= T(f_m)^\vee$, we obtain by Plancherel Theorem that $g_m \in L^2(G)$  
for each $m \in \N$. Moreover, from equation \eqref{76354433}, we have for any 
$m \in \mathbb{N}$
$$
\begin{aligned}
 \infty >\int_{G^\wedge} |f_m(\xi)|^2 d\nu(\xi)&= \int_{G^\wedge} (1 + \gamma(\xi)^2) \, |[T(f_m)](\xi)|^2 d\nu(\xi)
 \\[5pt]
  &= \int_{G^\wedge} (1 + \gamma(\xi)^2) \, |\widehat{g_m(\xi)}|^2 d\nu(\xi). 
\end{aligned}  
$$
Therefore, the sequence $\{g_m\}$ is uniformly bounded in $H^1_\gamma(G)$, with respect to $m \in \mathbb{N}$.  
By hypothesis there exists a subsequence of $\{g_m\}$, say $\{g_{m_j}\}$,
and a function $g \in L^2(G)$ such that, $g_{m_j}$ converges strongly to $g$ in $L^2(G)$ as $j \to \infty$. 
Consequently, we have 
$$T(f_{m_j})= \widehat{g_{m_j}} \to \widehat{g} \quad 
\text{in $L^2(G^\wedge)$ as $j \to \infty$},$$ that is, the operator $T$ is compact. 

\medskip
2. Now, let us assume that the operator $T$ is compact and then show that $H^1_\gamma(G) \subset \subset L^2(G)$. 
To this end, we take a sequence $\{f_m\}_{m\in \mathbb{N}}$ uniformly bounded in $H^1_\gamma(G)$. 
Then, due to the equivalence \eqref{87648764} there exists for each $m\in \mathbb{N}$, 
$g_m \in L^2(G^\wedge)$, such that $\widehat{f_m} = T(g_m)$. Thus for any $ m\in \mathbb{N}$, 
we have from equation \eqref{76354433} that
$$
\begin{aligned}
   \int_{G^\wedge} |g_m(\xi)|^2 d\nu(\xi)&= \int_{G^\wedge} (1 + \gamma(\xi)^2) \, |[T(g_m)](\xi)|^2 d\nu(\xi) 
   \\[5pt]
  & = \int_{G^\wedge} (1 + \gamma(\xi)^2) \, |\widehat{f_m(\xi)}|^2 d\nu(\xi)< \infty. 
\end{aligned}
$$
Then, the sequence $\{g_m\}$ is uniformly bounded in $L^2(G)$. Since the operator $T$ 
is compact, there exist $\{m_j\}_{j \in \mathbb{N}}$ and $g \in L^2(G^\wedge)$, such that  
$$
   \widehat{f_{m_j}}= T(g_{m_j})  \xrightarrow[j \to \infty]{} g \quad \text{in $L^2(G^\wedge)$}.
$$
Consequently, the subsequence $\{f_{m_j}\}$ converges to $g^\vee$ strongly in $L^2(G)$,
and thus $H^1_\gamma(G)$ is compactly embedded in $L^2(G)$.
\end{proof}

\begin{remark}
Due to Theorem \ref{876876872} the compactness characterization, that is
$H^1_\gamma(G) \subset \subset L^2(G)$, follows once
we show the conditions that the operator $T$ is compact. 
The study of the dual space of $G$, i.e. $G^\wedge$, and $\gamma$ it will be essential for this characterization.  
\end{remark}

\medskip
Recall from \eqref{CARACGCOMP} item $(ii)$ that, $G^\wedge$ is discrete since $G$ is compact. 
Then, $\nu$ is a countermeasure, and $\nu(\{\chi\})= 1$ for each singleton $\{\chi\}$, $\chi \in G^\wedge$. 
Now, for any $\chi \in G^\wedge$ fixed, we 
define the point mass function at $\chi$ by 
$$
   \delta_{\chi}(\xi):= 1_{\{\chi\}}(\xi),
\quad 
\text{for each $\xi \in G^\wedge$}. 
$$
Hence the set $\{\delta_\xi \; ; \; \xi \in G^\wedge \}$
is an orthonornal basis for $L^2(G^\wedge)$. Indeed, we first show the orthonormality. 
For each $\chi, \pi \in G^\wedge$, we have 
\begin{equation}
\label{87987948744}
    \langle \delta_\chi, \delta_\pi \rangle_{L^2(G^\wedge)}
     = \int_{G^\wedge} \delta_\chi(\xi)  \; \delta_\pi(\xi) \, d\nu(\xi)= \left\{
	\begin{array}{ccl}
		1, & \text{if} & \chi = \pi, 
		\\
		0, & \text{if} & \chi \not= \pi.
	\end{array}		 
\right.
\end{equation}
Now, let us show the density, that is $\overline{\{\delta_\xi \; ; \; \xi \in G^\wedge \}}= L^2(G^\wedge)$, or equivalently $\{\delta_\xi \; ; \; \xi \in G^\wedge \}^\perp= \{0\}$. 
For any $w \in \{\delta_\xi \; ; \; \xi \in G^\wedge \}^\perp$, we obtain 
$$
   0 =\langle \delta_\xi, w \rangle_{L^2(G^\wedge)} 
   = \int_{G^\wedge} \delta_\xi(\chi) w(\chi) \, d\nu(\chi) 
   = \int_{ \{ \xi \}} w(\chi) \, d\nu(\chi) = w(\xi)
$$
for any $\xi \in G^\wedge$, which proves the density. 

\medskip
From the above discussion, it is important to study the operator $T$
on elements of the set $\{\delta_\xi \; ; \; \xi \in G^\wedge \}$.
Then, we have the following 
\begin{theorem}
\label{876876876GG}
If the operator $T$ defined by \eqref{TCOMP} is compact,  then $G^\wedge$ is an enumerable set. 
\end{theorem}
\begin{proof} 1. First, let $\{\delta_\xi \; ; \; \xi \in G^\wedge \}$ be the orthonormal basis for $L^2(G^\wedge)$,
and $T$ the operator defined by \eqref{TCOMP}. Then, the function 
$\delta_\xi \in L^2(G^\wedge)$ is an eigenfunction of $T$ 
corresponding to the eigenvalue $(1+\gamma^2)^{-1/2}$, that is $\delta_\xi \neq 0$, and
\begin{equation}
\label{87486tydg}
    T(\delta_\xi)= \frac{\delta_\xi}{\sqrt{1+\gamma(\xi)^2}}.
\end{equation}

\medskip		
2. Now, since $T$ is compact and $\{\delta_\xi \; ; \; \xi \in G^\wedge \}$ is a basis for $L^2(G^\wedge)$, it must be enumerable from \eqref{87486tydg}. 
On the other hand, the function  $\xi \in G^\wedge \mapsto \delta_\xi \in L^2(G^\wedge)$ is injective, hence $G^\wedge$ is enumerable. 
\end{proof}

\begin{corollary}
If the operator $T$ defined by \eqref{TCOMP} is compact,  then
$L^2(G)$ is separable. 
\end{corollary}
	
\begin{proof} First, the Hilbert space $L^2(G^\wedge)$ is separable, since $\{\delta_\xi \; ; \; \xi\in G^\wedge\}$
is an enumerable orthonormal basis of it. Then, the proof follows applying the Plancherel Theorem. 
\end{proof}

\begin{corollary}
Let $G_B$ be the Bohr compactification of $\mathbb{R}^n$. 
Then $H^1_\gamma(G_B)$ is not compactly embedded in $L^2(G_B)$.
\end{corollary}
		
\begin{proof} Indeed, $G_B^\wedge$ is non enumerable.
\end{proof}
Consequently, $G^\wedge$ be enumerable is a necessarily condition for the operator $T$ be compact, which is not 
sufficient as shown by the Example \ref{NOSUFF} below. Indeed, it might depend on the chosen $\gamma$, see also  
Example \ref{NOSUFF10}. 

\medskip
To follow, we first recall the  
\begin{definition}
Let $G$ be a group (not necessarily a topological one) and $S$ a subset of it. 
The smallest subgroup of G containing every element of S, denoted $\langle S \rangle$, is called the subgroup 
generated by $S$. Equivalently, see Dummit, Foote \cite{DummitFoote} p.63, 
$$
  \langle S \rangle= \big\{ g^{\varepsilon_1}_1 g^{\varepsilon_2}_2 \ldots g^{\varepsilon_k}_k / 
  \text{$k \in \mathbb{N}$ and for each $j$, $g_j \in S, \varepsilon_j= \pm 1$} \big\}.
$$
Moreover, 
if a group $G= \langle S \rangle$, then $S$ is called a generator of $G$, and
in this case when S is finite, $G$ is called finitely generated. 
\end{definition} 	
	
\begin{theorem}
\label{876876876} 
If the operator $T$ defined by \eqref{TCOMP} is compact and
there exists a generator of $G^\wedge$ such that $\gamma$ is bounded on it, 
then $G^\wedge$ is finite generated.  
\end{theorem}
		
\begin{proof} 
Let $S_0$ be a generator of $G^\wedge$, such that $\gamma$ is bounded on it.  
Therefore, there exists $d_0 \geq 0$ such that,  
$$
   \text{for each $\xi \in S_0$, $\; \gamma(\xi) \leq d_0$}. 
$$
Now, since $T$ is compact and  $\Vert T \Vert \leq 1$, there exists  $0 < c \leq 1$ such that,
the set of eigenvectors 
$$
   \Big\{ \delta_\xi \; ; \; \xi \in G^\wedge \;\; \text{and} \;\; \frac{1}{\sqrt{1 + \gamma(\xi)^2}} \geq c \Big\}
   \equiv 
   \Big\{ \delta_\xi \; ; \; \xi \in G^\wedge \;\; \text{and} \;\; \gamma(\xi) \leq\sqrt{\frac{1}{c^2} - 1} \Big\}
$$
is finite, where we have used the Spectral Theorem for bounded compact operators. Therefore, since 
$$				
     \left\{ \delta_\xi \; ; \; \xi \in S_0 \right\} \subset 
     \left\{ \delta_\xi \; ; \; \xi \in G^\wedge \;\; \text{and} \;\; \gamma(\xi) \leq d_0 \right\}
$$
it follows that $S_0$ is a finite set, and thus $G^\wedge$ is finite generated. 
\end{proof}

\begin{example}[Infinite enumerable Torus] 
\label{NOSUFF}
Let us recall the Sobolev space $H^1_\gamma(\mathbb{T}^\N)$, where $\mathbb{T}^\N$ is the infinite enumerable Torus. 
We claim that: $H^1_\gamma(\mathbb{T}^\N)$ is not compactly embedded in $ L^2(\mathbb{T}^\N)$, 
for $\gamma$ defined in Exemple \ref{6576dgtftdefd}. 
Indeed, given $k \in \N$ we define $1_k \in \mathbb{Z}^\N$, such that it is zero for any coordinate 
$\ell \neq k$, and one in $k-$coordinate. Therefore, the set 
$$
    S_0 := \{ \xi_{1_k} \; ; \; k \in \N \}
$$
is an infinite generator of the dual group $(\mathbb{T}^\N)^\wedge$. 
Since for each $k \in \N$,  $\gamma(\xi_{1_k}) = 1$, i.e. bounded in $S_0$, applying Theorem \ref{876876876}
it follows that $H^1_\gamma(\mathbb{T}^\N)$ is not compactly embedded in $ L^2(\mathbb{T}^\N)$. 
\end{example}

\begin{remark} The above discussion in the Example \ref{NOSUFF} follows as well to the Sobolev space $H^1_\gamma(\mathbb{T}^I)$, where
$I$ is an index set (enumerable or not). Clearly, the Sobolev space $H^1_\gamma(\mathbb{T}^I)$ in not compactly embedded in $ L^2(\mathbb{T}^I)$, 
when $I$ is a non enumerable index set. Indeed, the set $(\mathbb{T}^I)^\wedge$ is non enumerable. 
\end{remark}

Now, we charactherize the condition on $\gamma: G^\wedge \to [0,\infty)$,
in order to $T$ be compact. More precisely, let us consider the following property:
\begin{equation}
\label{ConditionC}
 {\bf C}.  \quad \text{For each $d> 0$, the set 
$ \left\{ \xi \in G^\wedge \; ; \; \gamma(\xi) \leq d \right\}$
is finite}. 
\end{equation}
							
\begin{theorem}
\label{7864876874}
If $\gamma: G^\wedge \to [0,\infty)$ satisfies ${\bf C}$, then the operator $T$ defined by \eqref{TCOMP} is compact. 
\end{theorem}
		
\begin{proof}
By hypothesis, $\{ \xi \in G^\wedge \; ; \; \gamma(\xi) \leq d \}$ is finite, then we have
$$
    G^\wedge= \bigcup_{k \in \mathbb{N}} \left\{ \xi \in G^\wedge \; ; \; \gamma(\xi) \leq k \right\}. 
$$
Consequently, the set $G^\wedge$ is enumerable and we may write $G^\wedge= \{ \xi_i \}_{i \in \mathbb{N}}$. 

\medskip
Again, due to condition ${\bf C}$ for each $c \in (0,1)$ the set 
\begin{equation}
\label{868768767864120}
	\Big\{ \xi \in G^\wedge \; ; \; \frac{1}{\sqrt{1 + \gamma(\xi)^2}} \geq c \Big\}
\end{equation}
is finite. Since the function  $\xi \in G^\wedge \mapsto \delta_\xi \in L^2(G^\wedge)$ is injective,
the set $\{ \delta_{\xi_i} \; ; \; i\in \mathbb{N} \}$ 
is an enumerable orthonormal basis of eigenvectors for $T$, which corresponding eigenvalues satisfy 
$$
  \lim_{i \to \infty}  \frac{1}{\sqrt{1 + \gamma(\xi_i)^2}}= 0,
$$	
where we have used \eqref{868768767864120}. Consequently, $T$ is a compact operator. 
\end{proof}

\begin{example}[Bis: Infinite enumerable Torus]
\label{NOSUFF10}
There exists a function $\gamma_0$ such that, 
$H^1_{\gamma_0}(\mathbb{T}^\N)$ is compactly embedded in $ L^2(\mathbb{T}^\N)$. 
Indeed, we are going to show that, $\gamma_0$ satisfies ${\bf C}$. 
Let $\alpha \equiv (\alpha_\ell)_{\ell \in \mathbb{N}}$
be a sequence in $\R^\N$, such that for each $\ell \in \N$, $\alpha_\ell \geq 0$ and 
\begin{equation}
\label{weight}
\lim_{\ell \to \infty} \alpha_\ell = +\infty.
\end{equation}
Then, we define the following pseudo-metric in the dual group $(\mathbb{T}^\mathbb{N})^\wedge$ as follows
$$
   p_0(\xi_m, \xi_n):= 2\pi \sum_{\ell = 1}^\infty \alpha_\ell \;{\vert m_\ell - n_\ell \vert}, 
   \quad (m,n \in \mathbb{Z}^\mathbb{N}_{\rm c}),
$$
and consider $\gamma_0(\xi_m)= p_0(\xi_m,1)$. 
Thus for each $d> 0$, the set 
$$
    \{ m \in \mathbb{Z}^\mathbb{N}_{\rm c} \; ; \; \gamma_0(\xi_m) \leq d \} \quad \text{is finite.}
$$   
Indeed, from \eqref{weight} 
there exists $\ell_0 \in \N$, such that $\alpha_\ell> d$, for each $\ell \geq \ell_0$. 
Therefore, if $m \in \mathbb{Z}^\mathbb{N}_{\rm c}$ and the support of $m$ is not contained 
in $\{ 1, \ldots, \ell_0-1\}$, that is to say, there exists $\tilde{\ell} \geq \ell_0$, 
such that, $m_{\tilde{\ell}} \neq 0$. Then, 
$$
2\pi \sum_{\ell = 1}^\infty \alpha_\ell \;{\vert m_\ell  \vert}
\geq \alpha_{\tilde{\ell}} > d. 
$$
Consequently, we have 
$$
    \{ m \in \mathbb{Z}^\mathbb{N}_{\rm c} \; ; \; \gamma_0(\xi_m) \leq d \} 
    \subset 
    \{ m \in \mathbb{Z}^\mathbb{N}_{\rm c} \; ; \; {\rm supp} \ m \subset \{ 1, \ldots, \ell_0-1\} \},
$$   
which is a finite set. Finally, applying Theorem \ref{7864876874} we obtain that,
the Sobolev space 
$H^1_{\gamma_0}(\mathbb{T}^\N)$ is compactly embedded in $ L^2(\mathbb{T}^\N)$. 
\end{example}

\subsection{Application. On a class of quasiperiodic functions}
\label{4563tgf5fd3}

In this section we consider the important class 
of quasiperiodic functions (see \cite{JikovKozlovOleinik}), which 
includes the class of $[0,1)^n-$periodic functions. 
 
 \smallskip
Let $\lambda_1,\lambda_2,\ldots,\lambda_m \in \mathbb{R}^n$ be $m-$linear independent 
vectors with respect to $\mathbb{Z}$, and consider the following matrix
\begin{equation*}
\Lambda := {\left(
\begin{array}{c}
	\lambda_1
\\
	\lambda_2
\\
	\vdots 
\\
	\lambda_m
\end{array}
\right)}_{m\times n}
\end{equation*}
such that, for each $d> 0$ the set 
\begin{equation}
\label{7863948tyfedf}
   \{ k \in \mathbb{Z}^m \; ; \; {\vert \Lambda^T k \vert} \leqslant d\} \quad \text{is finite.}
\end{equation}
Therefore, we are considering the class of quasiperiodic functions satisfying 
condition \eqref{7863948tyfedf}. This set is not empty, for instance let us define 
the matrix $B:= \Lambda \Lambda^T$, such that $\det B> 0$, which is called here
positive quasi-periodic functions. It is not difficult to see that, positive quasiperiodic functions 
satisfies \eqref{7863948tyfedf}. 
Indeed, it is sufficiently to observe that, for each $k \in \mathbb{Z}^m$, we have 
$$
   |k|= | B^{-1} B k | \leq \|B^{-1}\| \|\Lambda\| |\Lambda^T k|. 
$$
Moreover, since $\lambda_1,\lambda_2,\ldots,\lambda_m \in \mathbb{R}^n$ are $m-$linear independent 
vectors with respect to $\mathbb{Z}$, (this property does not imply $\det B> 0$), 
the dynamical system $\tau: \mathbb{R}^n \times \mathbb{T}^m \to \mathbb{T}^m$,  given by 
\begin{equation}
\label{6973846tyd4f54e54}
   \tau(x)\omega := \omega + \Lambda x - \left\lfloor \omega + \Lambda x \right\rfloor
\end{equation}
is ergodic. 

\medskip
Now we remark that, the application 
${ \varphi : \mathbb{R}^n \to \mathbb{T}^m }$, $ \varphi(x) := \Lambda x - \left\lfloor \Lambda x \right\rfloor$,
is a continuous homeomorphism of groups. Then, we have
$$
   \tau(x)\omega = \varphi(x) \omega \equiv \omega + \Lambda x - \left\lfloor \omega + \Lambda x \right\rfloor. 
$$
Consequently, under the conditions of the previous sections, we obtain for each 
$ k\in \mathbb{Z}^m$, $\gamma(\xi_k)= 2\pi {\vert \Lambda^T k \vert}$, 
and applying Theorem \ref{7864876874} (recall \eqref{7863948tyfedf}),
it follows that 
\begin{equation*}
	H^1_\gamma {\left( \mathbb{T}^m \right)} \subset \! \subset L^2{\left( \mathbb{T}^m \right)}.
\end{equation*}
Therefore, given a stochastic deformation $\Phi$, we have $\mathcal{H}_\Phi \subset \! \subset \mathcal{L}_\Phi$
for the class of quasiperiodic functions satisfying \eqref{7863948tyfedf}, 
and it follows a solution of the Bloch's spectral cell equation.  

\section*{Acknowledgements}
Conflict of Interest: Author Wladimir Neves has received research grants from CNPq
through the grant  308064/2019-4, and also by FAPERJ 
(Cientista do Nosso Estado) through the grant E-26/201.139/2021. 
Author Jean Silva has received research grants from CNPq through the Grant 303533/2020-0.

\end{document}